\documentclass{amsart}

\begin{document}
\newtheorem{theorem}{Theorem}[section]
\newtheorem{lemma}{Lemma}[section]
\newtheorem{prop}{Proposition}[section]
\newtheorem{corollary}{Corollary}[section]
\newtheorem{example}{Example}[section]
\newtheorem{definition}{Definition}[section]

\numberwithin{equation}{section}

\title{Partition's sensitivity for measurable maps}

\author[C. A. Morales]{C. A. Morales}
\address{Instituto de Mat\'ematica, UFRJ,
P. O. Box 68530, 21945-970 Rio de Janeiro, Brazil.}
\email{morales@impa.br}

\thanks{Partially supported by CNPq, CAPES-Prodoc, FAPERJ and PRONEX/DS from Brazil.}


\begin{abstract}
We study countable partitions for measurable maps on measure spaces such that for all point $x$
the set of points with the same itinerary of $x$ is negligible.
We prove that in nonatomic probability spaces every strong generator (Parry, W.,
{\em Aperiodic transformations and generators},
J. London Math. Soc. 43 (1968), 191--194)
satisfies this property but not conversely.
In addition, measurable maps carrying partitions with this property are
aperiodic and their corresponding spaces are nonatomic.
From this we obtain a characterization of nonsingular countable to one mappings with these partitions
on nonatomic Lebesgue probability spaces
as those having strong generators.
Furthermore, maps carrying these partitions
include the ergodic measure-preserving ones with positive entropy on probability spaces
(thus extending a result in Cadre, B., Jacob, P.,
{\em On pairwise sensitivity},
J. Math. Anal. Appl. 309 (2005), no. 1, 375--382).
Some applications are given.
\end{abstract}

\maketitle

\section{Introduction}\label{sec1}

\noindent
In this paper we study countable partitions $P$ for measurable maps $f: X\to X$
on measure spaces $(X,\mathcal{B},\mu)$ such that
for all $x\in X$ {\em the set of points with the same itinerary of $x$ is negligible}.
In other words,
\begin{equation}
\label{def-meas-sen}
\mu(\{y\in X:f^n(y)\in P(f^n(x)),\quad \forall n\in\mathbb{N}\})=0,
\quad\forall x\in X,
\end{equation}
where $P(x)$ stands for the element of $P$ containing $x\in X$.
For simplicity, we call these partitions {\em measure-sensitive partitions}.

We prove that in a nonatomic probability space every {\em strong generator}
is a measure-sensitive partition but not conversely
(results about strong generators can be found in \cite{hs}, \cite{jk}, \cite{k}, \cite{p}, \cite{p2} and \cite{p11}).
We also exhibit examples of measurable maps
in nonatomic probability spaces carrying measure-sensitive partitions which are not strong generators.
Motivated by these examples we shall study measurable maps on measure spaces
carrying measure-sensitive partitions (called {\em measure-expansive maps} for short).
Indeed, we prove that every measure-expansive map is aperiodic and also, in the probabilistic case, that
its corresponding space is nonatomic.
From this we obtain a characterization of nonsingular countable to one measure-expansive mappings
on nonatomic Lebesgue probability spaces
as those having strong generators.
Furthermore, we prove that every ergodic measure-preserving map with positive entropy is a probability space
is measure-expansive (thus extending a result in \cite{cj}).
As an application we obtain some properties
for ergodic measure-preserving maps with positive entropy
(c.f. corollaries \ref{nonatomic} and \ref{eventually-aperiodic}).

\section{Statements and proofs}

\noindent
Hereafter the term {\em countable} will mean either finite or countably infinite.

A {\em measure space} is a triple
$(X,\mathcal{B},\mu)$ where $X$ is a set, $\mathcal{B}$ is a $\sigma$-algebra of subsets of
$X$ and $\mu$ is a positive measure in $\mathcal{B}$.
A {\em probability space} is one for which $\mu(X)=1$.
A measure space is {\em nonatomic} if it has no
{\em atoms}, i.e., measurable sets $A$ of positive measure
satisfying $\mu(B)\in\{0,\mu(A)\}$ for every measurable set $B\subset A$.
A {\em partition} is
a disjoint collection $P$ of nonempty measurable sets
whose union is $X$. We allow $\mu(\xi)=0$ for some $\xi\in P$.
If $f: X\to X$ is measurable and $k\in \mathbb{N}$ we define
$f^{-k}(P)=\{f^{-k}(\xi):\xi\in P\}$ which is a (countable) partition if $P$ is.
A {\em strong generator} of $f$
is a countable partition $P$ for which the smallest
$\sigma$-algebra of $\mathcal{B}$ containing $\bigcup_{k\in\mathbb{N}}f^{-k}(P)$ equals $\mathcal{B}$
(mod $0$) (see \cite{p}).

The result below is the central motivation of this paper.

\begin{theorem}
\label{strong-gen}
Every strong generator of a measurable map in a nonatomic probability space
is a measure-sensitive partition.
\end{theorem}

\proof
Recall that the join of finitely many partitions $P_0,\cdots,P_n$
is the partition defined by
$$
\bigvee_{k=0}^nP_k=\left\{\bigcap_{k=0}^n\xi_k:\xi_k\in P_k, \forall 0\leq k\leq n\right\}.
$$
Given partitions $P$ and $Q$ we write
$P\leq Q$ to mean that each member of $Q$ is contained in some member of $P$ (mod $0$).
A sequence of partitions $\{P_n:n\in\mathbb{N}\}$ (or simply $P_n$) is {\em increasing}
if $P_0\leq P_1\leq \cdots\leq P_n\leq\cdots$.
Certainly
\begin{equation}
\label{obama1}
P_n=\bigvee_{k=0}^nf^{-k}(P),
\quad\quad n\in\mathbb{N},
\end{equation}
defines an increasing sequence of countable partitions
satisfying
$$
P_n(x)=\bigcap_{k=0}^nf^{-k}(P(f^k(x)),
\quad\quad\forall x\in X.
$$
Since for all $x\in X$ one has
$$
\{y\in X:f^n(y)\in P(f^n(x)),\quad\forall n\in\mathbb{N}\}=\bigcap_{n=0}^\infty f^{-n}(P(f^n(x)))
=\bigcap_{n=0}^\infty P_n(x),
$$
we obtain that the identity below
\begin{equation}
 \label{*}
\lim_{n\to\infty}\sup_{\xi\in P_n}\mu(\xi)=0
\end{equation}
implies (\ref{def-meas-sen}).

Now suppose that $P$ is a strong generator
of a measurable map $f: X\to X$ in a nonatomic probability space
$(X,\mathcal{B},\mu)$.
Then, the sequence (\ref{obama1})
generates $\mathcal{B}$ (mod $0$).
From this and Lemma 5.2 p. 8 in \cite{man} we obtain that
the set of all finite unions of
elements of these partitions is everywhere dense in the measure algebra
associated to $(X,\mathcal{B},\mu)$.
Consequently, Lemma 9.3.3 p. 278 in \cite{b} implies that the sequence (\ref{obama1}) satisfies (\ref{*})
and then (\ref{def-meas-sen}) holds.
\endproof

We shall see later in Example \ref{measure-strong} that
the converse of this theorem is false,
i.e., there are certain measurable maps in nonatomic probability spaces
carrying measure-sensitive partitions which are not strong generators.
These examples motivates the study of measure-sensitive partitions for measurable maps
in measure spaces. For this we use the following auxiliary
concept motivated by the notion of Lebesgue sequence of partitions
(c.f. p. 81 in \cite{man}).

\begin{definition}
 \label{ms-partition}
A {\em measure-sensitive sequence of partitions} of $(X,\mathcal{B},\mu)$ is
an increasing sequence of countable partitions $P_n$
such that $\mu\left(\bigcap_{n\in\mathbb{N}}\xi_n\right)=0$ for all sequence
of measurable sets $\xi_n$ satisfying $\xi_n\in P_n$, $\forall n\in\mathbb{N}$.
A {\em measure-sensitive space} is a measure space carrying measure-sensitive sequences of partitions.
\end{definition}

At first glance we observe that (\ref{*}) is sufficient
condition for an increasing sequence $P_n$ of countable partitions
to be measure-sensitive (it is also necessary in probability spaces).
On the other hand, the class of measure-sensitive spaces is broad enough to include all
nonatomic standard probability spaces. Recall that a standard probability space is a probability space $(X,\mathcal{B},\mu)$
whose underlying measurable space $(X,\mathcal{B})$ is isomorphic to a Polish
space equipped with its Borel $\sigma$-algebra (e.g. \cite{aa}).
Precisely we have the following proposition.

\begin{prop}
\label{--}
All nonatomic standard probability spaces are measure-sensitive.
\end{prop}

\proof
As is well-known, for every nonatomic standard probability space $(X,\mathcal{B},\mu)$
there are a measurable subset
$X_0\subset X$ with $\mu(X\setminus X_0)=0$ and a sequence of countable partitions $Q_n$ of $X_0$
such that $\bigcap_{n\in\mathbb{N}}\xi_n$ contains at most
one point for every sequence of measurable sets $\zeta_n$ in $X_0$
satisfying $\zeta_n\in Q_n$, $\forall n\in\mathbb{N}$ (c.f. \cite{man} p. 81).
Defining $P_n=\{X\setminus X_0\}\cup Q_n$
we obtain an increasing sequence of countable partitions of $(X,\mathcal{B},\mu)$.
It suffices to prove that this sequence is measure-sensitive.
For this take a fixed (but arbitrary)
sequence of measurable sets $\xi_n$ of $X$ with $\xi_n\in P_n$ for all $n\in\mathbb{N}$.
It follows from the definition of $P_n$ that
either $\xi_n=X\setminus X_0$ for some
$n\in \mathbb{N}$, or,
$\xi_n\in Q_n$ for all $n\in\mathbb{N}$.
Then, the intersection
$\bigcap_{n\in\mathbb{N}}\xi_n$ either is contained in $X\setminus X_0$
or reduces to a single measurable point.
Since both $X\setminus X_0$ and the measurable points have measure zero
(for nonatomic spaces are diffuse \cite{b}) we obtain
$\mu\left(\bigcap_{n\in\mathbb{N}}\xi_n\right)=0$.
As $\xi_n$ is arbitrary we are done.
\endproof

Although measure-sensitive probability spaces need not be standard we have that all of them are
nonatomic. Indeed, we have the following result of later usage.

\begin{prop}
 \label{ms-nonatomic}
All measure-sensitive probability spaces are nonatomic.
\end{prop}

\proof
Suppose by contradiction that there is a measure-sensitive probability space
$(X,\mathcal{B},\mu)$ with an atom $A$. Take a measure-sensitive sequence of partitions $P_n$.
Since $A$ is an atom one has that $\forall n\in\mathbb{N}$
$\exists !\xi_n\in P_n$ such that $\mu(A\cap \xi_n)>0$
(and so $\mu(A\cap \xi_n)=\mu(A)$).
Notice that $\mu(\xi_n\cap \xi_{n+1})>0$ for, otherwise,
$\mu(A)\geq \mu(A\cap (\xi_n\cup \xi_{n+1}))=\mu(A\cap \xi_n)+\mu(A\cap \xi_{n+1})=2\mu(A)$
which is impossible in probability spaces.
Now observe that $\xi_n\in P_n$ and $P_n\leq P_{n+1}$, so, there is $L\subset P_{n+1}$ such that
\begin{equation}
\label{zeta}
 \mu\left(\xi_n\bigtriangleup \bigcup_{\zeta\in L}\zeta\right)=0.
\end{equation}
If $\xi_{n+1}\cap \left(\bigcup_{\zeta\in L}\zeta\right)=\emptyset$
we would have
$\xi_n\cap \xi_{n+1}=\xi_n\cap\xi_{n+1}\setminus\bigcup_{\zeta\in L}\zeta$
yielding
$$
\mu(\xi_n\cap \xi_{n+1})=\mu\left(\xi_n\cap\xi_{n+1}\setminus\bigcup_{\zeta\in L}\zeta\right)
\leq \mu\left(\xi_n\setminus \bigcup_{\zeta\in L}\zeta\right)=0
$$
which is absurd.
Hence $\xi_{n+1}\cap \left(\bigcup_{\zeta\in L}\zeta\right)\neq\emptyset$
and then $\xi_{n+1}\in L$ for $P_{n+1}$ is a partition and $\xi_{n+1}\in P_{n+1}$.
Using (\ref{zeta}) we obtain $\xi_{n+1}\subset \xi_n$ (mod $0$) so
$A\cap \xi_{n+1}\subset A\cap \xi_{n}$ (mod $0$) for all $n\in\mathbb{N}^+$.
From this and well-known properties of probability spaces we obtain
$$
\mu\left(A\cap \bigcap_{n\in\mathbb{N}}\xi_n\right)
=
\mu\left(\bigcap_{n\in\mathbb{N}}(A\cap \xi_n)\right)
=\lim_{n\to\infty}\mu(A\cap \xi_n)=\mu(A)>0.
$$
But $P_n$ is measure-sensitive
and $\xi_n\in P_n$, $\forall n\in \mathbb{N}$, so
$\mu\left(\bigcap_{n\in\mathbb{N}}\xi_n\right)=0$
yielding $\mu\left(A\cap \bigcap_{n\in\mathbb{N}}\xi_n\right)=0$
which contradicts the above expression.
This contradiction yields the proof.
\endproof

The following equivalence relates both measure-sensitive partitions and
measure-sensitive sequences of partitions.

\begin{lemma}
\label{l0}
The following properties are equivalent for measurable maps
$f: X\to X$ and countable partitions $P$ on measure spaces $(X,\mathcal{B},\mu)$:
\begin{enumerate}
 \item[(i)]
The sequence $P_n$ in (\ref{obama1}) is measure-sensitive.
\item[(ii)]
The partition $P$ is measure-sensitive.
\item[(iii)]
The partition $P$ satisfies
$$
\mu(\{y\in X:f^n(y)\in P(f^n(x)),\quad \forall n\in \mathbb{N}\})=0,
\quad \forall \mu\mbox{-a.e. }x\in X.
$$
\end{enumerate}
\end{lemma}

\proof
Previously we state some notation.

Given a partition $P$ and $f: X\to X$ measurable
we define
$$
P_\infty(x)=\{y\in X:f^n(y)\in P(f^n(x)),\forall n\in\mathbb{N}\},
\quad\quad \forall x\in X.
$$
Notice that
\begin{equation}
\label{pinfinity}
P_\infty(x)=\bigcap_{n\in\mathbb{N}^+}P_n(x)
\end{equation}
and
\begin{equation}
\label{pn}
P_n(x)=\bigcap_{i=0}^{n}f^{-i}(P(f^i(x)))
\end{equation}
so each $P_\infty(x)$ is a measurable set.
For later use we keep the following identity
\begin{equation}
\label{lee1}
\left(\bigvee_{i=0}^nf^{-i}(P)\right)(x)=P_n(x),
\quad\quad\forall x\in X.
\end{equation}
Clearly (\ref{def-meas-sen}) (resp. (iii)) is equivalent to
$\mu(P_\infty(x))=0$ for every $x\in X$ (resp. for $\mu$-a.e. $x\in X$).

First we prove that (i) implies (ii).
Suppose that the sequence (\ref{obama1}) is measure-sensitive and
fix $x\in X$.
By (\ref{pinfinity}) and (\ref{lee1}) we have
$P_\infty(x)=\bigcap_{n\in\mathbb{N}}\xi_n$
where $\xi_n=P_n(x)\in P_n$. As the sequence $P_n$ is measure-sensitive we
obtain $\mu(P_\infty(x))=\mu\left(\bigcap_{n\in\mathbb{N}}\xi_n\right)=0$ proving (ii).
Conversely, suppose that (ii) holds and let $\xi_n$ be a sequence of measurable sets with
$\xi_n\in P_n$ for all $n$.
Take $y\in \bigcap_{n\in\mathbb{N}}\xi_n$.
It follows that
$y\in P_n(x)$ for all $n\in\mathbb{N}$ whence
$y\in P_\infty(x)$ by (\ref{obama1}).
We conclude that
$\bigcap_{n\in\mathbb{N}}\xi_n\subset P_\infty(x)$ therefore
$\mu\left(\bigcap_{n\in\mathbb{N}}\xi_n\right)\leq \mu(P_\infty(x))=0$
proving (i).

To prove that (ii) and (iii) are equivalent
we only have to prove that (iii) implies (i).
Assume by contradiction that $P$ satifies (iii) but not (ii).
Since $\mu$ is a probability and (3) holds the set
$X'=\{x\in X:\mu(P_\infty(x))=0\}$ has measure one.
Since (ii) does not hold there is $x\in X$ such that $\mu(P_\infty(x))>0$.
Since $\mu$ is a probability and $X'$ has measure one we would
have $P_\infty(x)\cap X'\neq\emptyset$ so
there is $y\in P_\infty(x)$ such that $\mu(P_\infty(y))=0$.
But clearly the collection $\{P_\infty(x):x\in X\}$
is a partition (for $P$ is) so
$P_\infty(x)=P_\infty(y)$ whence $\mu(P_\infty(x))=\mu(P_\infty(y))=0$ which is a contradiction.
This ends the proof.
\endproof

Recall that a measurable map $f: X\to X$ is
{\em measure-preserving} if $\mu\circ f^{-1}=\mu$. Moreover, it is
{\em ergodic} if every measurable invariant set $A$ (i.e. $A=f^{-1}(A)$ (mod $0$))
satisfies either $\mu(A)=0$ or $\mu(X\setminus A)=0$; and
{\em totally ergodic} if $f^n$ is ergodic for all $n\in \mathbb{N}^+$.

\begin{example}
\label{re2'}
If $f$ is a totally ergodic measure-preserving map of a probability space,
then every countable partition $P$
with $0<\mu(\xi)<1$ for some $\xi\in P$ is measure-sensitive with respect to
$f$ (this follows from the equivalence (iii) in Lemma \ref{l0} and Lemma 1.1 p. 208 in \cite{man}).
\end{example}

Hereafter we fix a measure space $(X,\mathcal{B},\mu)$ and a measurable map $f: X\to X$.
We shall not assume that $f$ is measure-preserving unless otherwise stated.

Using the Kolmogorov-Sinai's entropy we obtain
sufficient conditions for the measure-sensitivity of a given partition.
Recall that the {\em entropy} of a finite partition $P$ is defined by
$$
H(P)=-\sum_{\xi\in P}\mu(\xi)\log\mu(\xi).
$$
The {\em entropy} of a finite partition $P$ with respect to
a measure-preserving map $f$ is defined by
$$
h(f,P)=\lim_{n\to\infty}\frac{1}{n}
H(P_{n-1}).
$$
Then, we have the following lemma.

\begin{lemma}
 \label{shannon}
A finite partition with finite positive entropy of
an ergodic measure-preserving map in a probability space is measure-sensitive.
\end{lemma}

\proof
Since the map $f$ (say) is ergodic, the Shannon-McMillan-Breiman Theorem (c.f. \cite{man} p. 209) implies
that the partition $P$ (say) satisfies
\begin{equation}
 \label{smb}
-\lim_{n\to\infty}\frac{1}{n}\log(\mu(P_n(x)))=h(f,P),
\quad\quad\mu\mbox{-a.e. }x\in X,
\end{equation}
where $P_n(x)$ is as in (\ref{pn}).
On the other hand, $P_{n+1}(x)\subset P_n(x)$ for all $n$
so (\ref{pinfinity}) implies
\begin{equation}
 \label{lll}
\mu(P_\infty(x))=\lim_{n\to\infty}\mu(P_n(x)),
\quad\quad\forall x\in X.
\end{equation}
But
$h(f,P)>0$
so (\ref{smb}) implies that $\mu(P_n(x))$ goes to zero
for $\mu$-a.e. $x\in X$.
This together with (\ref{lll}) implies
that $P$ satisfy the equivalence (iii) in Lemma \ref{l0} so $P$ is measure-sensitive.
\endproof

In the sequel we study measurable maps
carrying measure-sensitive partitions
(we call them {\em measure-expansive maps} for short).
It follows at once from
Lemma \ref{l0} that these maps only exist on
measure-sensitive spaces.
Consequently we obtain the following result from
Proposition \ref{ms-nonatomic}.

\begin{theorem}
\label{theorem1}
A probability space carrying measure-expansive maps
is nonatomic.
\end{theorem}

A simple but useful example is as follows.

\begin{example}
 \label{re2}
The irrational rotations in the circle are measure-expansive maps with respect to the Lebesgue measure.
This follows from Example \ref{re2'} since all such maps are measure-preserving and totally ergodic.
\end{example}

On the other hand, it is not difficult to find examples of
measure-expansive measure-preserving maps which are not ergodic.
These examples
together with Example \ref{re2}
suggest the question whether an ergodic measure-preserving map is measure-expansive.
However, the answer is negative by the following example.

\begin{example}
\label{flu}
A measure-sensitive partition has necessarily more than one element.
Consequently, if $\mathcal{B}=\{X,\emptyset\}$ then no map is measure-expansive
although they are all ergodic measure-preserving.
\end{example}

Despite of this we still can give conditions for the measure-expansivity of ergodic
measure-preserving maps as follows.

Recall that the {\em entropy} (c.f. \cite{man}, \cite{w}) of $f$ is defined by
$$
h(f)=\sup\{h(f,Q):Q\mbox{ is a finite partition of }X\}.
$$
We obtain a natural generalization of Theorem 3.1 in \cite{cj}.

\begin{theorem}
 \label{cad-jac}
Ergodic measure-preserving maps with positive entropy
in probability spaces are measure-expansive.
\end{theorem}

\proof
Let $f$ be one of such a map with entropy $h(f)>0$.
We can assume that $h(f)<\infty$.
It follows that there is a finite partition $Q$ with $0<h(f,Q)<\infty$.
Taking
$P=\bigvee_{i=0}^{n-1}f^{-i}(Q)$ with $n$ large we obtain a finite partition with finite positive entropy
since $h(f,P)=h(f,Q)>0$.
It follows that $P$ is measure-sensitive by Lemma \ref{shannon} whence
$f$ is measure-expansive by definition.
\endproof

A first consequence of the above result
is that the converse of Theorem \ref{strong-gen} is false.

\begin{example}
\label{measure-strong}
Let $f: X\to X$ be a homeomorphism with positive topological entropy
of a compact metric space $X$.
By the variational principle \cite{w} there is a Borel probability measures $\mu$ with respect to which
$f$ is an ergodic measure-preserving map with positive entropy.
Then, by Theorem \ref{cad-jac}, $f$ carries a measure-sensitive partition
which, by Corollary 4.18.1 in \cite{w}, cannot be a strong generator.
Consequently, there are measurable maps
in certain nonatomic probability spaces carrying measure-sensitive partitions which are not strong generators.
\end{example}

On the other hand, it is also false that ergodic measure-expansive measure-preserving maps
on probability spaces have positive entropy.
The counterexamples are precisely the irrational circle rotations (c.f. Example \ref{re2}).
Theorems \ref{theorem1} and \ref{cad-jac} imply the probably well-known result below.

\begin{corollary}
\label{nonatomic}
Probability spaces carrying ergodic measure-preserving maps
with positive entropy are nonatomic.
\end{corollary}

In the sequel we analyse the aperiodicity of measure-expansive maps.
According to \cite{p} a measurable map $f$ is {\em aperiodic} whenever for all $n\in\mathbb{N}^+$
if $n\in\mathbb{N}^+$ and $f^n(x)=x$ on a measurable set $A$, then $\mu(A)=0$.
Let us extend this definition in the following way.

\begin{definition}
\label{ea}
We say that $f$ is {\em eventually aperiodic} whenever the following property holds for every
$(n,k)\in \mathbb{N}^+\times \mathbb{N}$:
If $A$ is a measurable set such that
for every $x\in A$ there is $0\leq i\leq k$
such that $f^{n+i}(x)=f^i(x)$, then $\mu(A)=0$.
\end{definition}

It follows easily from the definition that an eventually periodic map is aperiodic.
The converse is true for invertible maps but not in general
(e.g. the constant map $f(x)=c$ where $c$ is a measurable point of zero mass).

With this definition we can state the following result.

\begin{theorem}
\label{p1}
Every measure-expansive map is eventually aperiodic
(and so aperiodic).
\end{theorem}

\proof
Let $f$ be a measure-expansive map of $X$.
Take $(n,k)\in \mathbb{N}^+\times \mathbb{N}$ and a measurable set $A$
such that for every $x\in A$ there is $0\leq i\leq k$ such that
$f^{n+i}(x)=f^i(x)$.
Then,
\begin{equation}
\label{equino1}
A\subset \bigcup_{i=0}^kf^{-i}(Fix(f^n)),
\end{equation}
where
$Fix(g)=\{x\in X:g(x)=x\}$ denotes the set of fixed points of a map $g$.

Let $P$ be a measure-sensitive partition of $f$.
Then, $\bigvee_{m=0}^{k+n}f^{-m}(P)$ is a countable partition.
Fix $x,y\in A\cap \xi$. In particular $\xi=\left(\bigvee_{m=0}^{k+n}f^{-m}(P)\right)(x)$
whence $y\in \left(\bigvee_{m=0}^{k+n}f^{-m}(P)\right)(x)$.
This together with (\ref{pn}) and (\ref{lee1})
yields
\begin{equation}
\label{la}
f^m(y)\in P(f^m(x)),
\quad\quad\forall 0\leq m\leq k+n.
\end{equation}
But $x,y\in A$ so (\ref{equino1}) implies
$f^i(x),f^j(y)\in Fix(f^n)$ for some $i,j\in\{0,\cdots, k\}$.
We can assume that $j\geq i$ (otherwise we interchange the roles
of $x$ and $y$ in the argument below).

Now take $m>k+n$. Then, $m>j+n$ so $m-j=pn+r$ for some $p\in\mathbb{N}^+$ and some integer $0\leq r<n$.
Since $0\leq j+r<k+n$ (for $0\leq j\leq k$ and $0\leq r<n$) one gets
\begin{eqnarray*}
f^m(y)=f^{m-j}(f^j(y)) & = & f^{pn+r}(f^j(y)) \\
& = & f^r(f^{pn}(f^j(y))) \\
& = &  f^{j+r}(y)\\
& \stackrel {(\ref{la})}{\in} & P(f^{j+r}(x)).
\end{eqnarray*}
But
\begin{eqnarray*}
P(f^{j+r}(x))= P(f^{j+r-i}(f^i(x))) & = & P(f^{j+r-i}(f^{pn}(f^i(x)))) \\
& = &  P(f^{m-i}(f^i(x)))\\
& = & P(f^m(x))
\end{eqnarray*}
so
$$
f^m(y)\in P(f^m(x)),
\quad\quad\forall m>k+n.
$$
This together with (\ref{la}) implies that $f^m(y)\in P(f^m(x))$ for all $m\in\mathbb{N}$ whence
$y\in P_\infty(x)$.
Consequently
$
A\cap \xi\subset P_\infty(x).
$
As $P$ is measure-sensitive Lemma \ref{l0} implies
$$
\mu(A\cap\xi)=0,
\quad\quad\forall \xi\in \bigvee_{i=0}^{k+n}f^{-i}(P).
$$
On the other hand, $\bigvee_{i=0}^{k+n}f^{-i}(P)$ is a partition so
$$
A=\bigcup_{\xi\in \bigvee_{i=0}^{k+n}f^{-i}(P)}(A\cap \xi)
$$
and then
$\mu(A)=0$ since $\bigvee_{i=0}^{k+n}f^{-i}(P)$ is countable.
This ends the proof.
\endproof

It follows from Lemma \ref{strong-gen} that, in nonatomic probability spaces, every measurable map
carrying strong generators is measure-expansive.
This motivates the question as to whether every measure-expansive map has a strong generator.
We give a partial positive answer for certain maps defined as follows.
We say that $f$ is {\em countable to one (mod $0$)} if
$f^{-1}(x)$ is countable for $\mu$-a.e. $x\in X$.
We say that $f$ is {\em nonsingular} if a measurable set $A$ has measure zero
if and only if $f^{-1}(A)$ also does.
All measure-preserving maps are nonsingular.
A {\em Lebesgue probability space} is a complete measure space which is isomorphic
to the completion of a standard probability space (c.f. \cite{aa}, \cite{b}).

\begin{corollary}
 \label{th1}
The following properties are equivalent
for nonsingular countable to one (mod $0$) maps $f$ on nonatomic Lebesgue probability spaces:

\begin{enumerate}
 \item
$f$ is measure-expansive.
\item
$f$ is eventually aperiodic.
\item
$f$ is aperiodic.
\item
$f$ has a strong generator.
\end{enumerate}
\end{corollary}

\proof
Notice that (1) $\Rightarrow$ (2) by Theorem \ref{p1} and (2) $\Rightarrow$ (3)
follows from the definitions.
On the other hand, (3) $\Rightarrow$ (4)
by a Parry's theorem (c.f. \cite{p}, \cite{p11}, \cite{p2}) while (4) $\Rightarrow$ (1)
by Lemma \ref{strong-gen}.
\endproof

Denote by $Fix(g)=\{x\in X:g(x)=x\}$ the set of fixed points of a mapping $g$.

\begin{corollary}
 \label{ccc1}
If $f^k=f$ for some integer $k\geq 2$, then $f$ is not measure-expansive.
\end{corollary}

\proof
Suppose by contradiction that it does.
Then, $f$ is eventually aperiodic by Theorem \ref{p1}.
On the other hand,
if $x\in X$ then $f^k(x)=f(x)$ so
$f^{k-1}(f^k(x))=f^{k-1}(f(x))=f^k(x)$ therefore
$f^k(x)\in Fix(f^{k-1})$
whence $X\subset f^{-k}(Fix(f^{k-1}))$.
But since $f$ is eventually aperiodic, $n=k-1\in\mathbb{N}^+$
and $X$ measurable we obtain from the definition that
$\mu(X)=0$ which is absurd.
This ends the proof.
\endproof

\begin{example}
\label{esquisito1}
By Corollary \ref{ccc1} neither the identity $f(x)=x$ nor the constant map $f(x)=c$
are measure-expansive (for they satisfy $f^2=f$).
In particular, the converse of Theorem \ref{p1} is false
for the constant maps are eventually aperiodic but not measure-expansive.
\end{example}

It is not difficult to prove that an ergodic measure-preserving map of a nonatomic probability space
is aperiodic.
Then, Corollary \ref{nonatomic} implies
the well-known fact that {\em all ergodic measure-preserving maps with positive entropy
on probability spaces are aperiodic}.
However, using theorems \ref{cad-jac} and \ref{p1} we obtain
the following stronger result.

\begin{corollary}
\label{eventually-aperiodic}
All ergodic measure-preserving maps with positive entropy on probability spaces
are {\em eventually aperiodic}.
\end{corollary}

Now we study the following variant of aperiodicity introduced
in \cite{hs} p. 180.

\begin{definition}
We say that $f$ is {\em aperiodic*} (\footnote{called aperiodic in \cite{hs}.}) whenever
for every measurable set of positive measure $A$ and $n\in\mathbb{N}^+$ there
is a measurable subset $B\subset A$ such that $\mu(B\setminus f^{-n}(B))>0$.
\end{definition}

Notice that aperiodicity* implies the aperiodicity used in \cite{jk} or \cite{s}
(for further comparisons see p. 88 in \cite{k}).

On the other hand, a measurable map $f$ is {\em negative nonsingular} if
$\mu(f^{-1}(A))=0$ whenever $A$ is a measurable set with $\mu(A)=0$.
Some consequences of the aperiodicity* on negative nonsingular maps
in probability spaces are given in \cite{k}.
Observe that every measure-preserving map is negatively nonsingular.

Let us present two technical (but simple) results for later usage.
We call a measurable set $A$ satisfying
$A\subset f^{-1}(A)$ (mod $0$) a
{\em positively invariant set (mod $0$)}.
For completeness we prove the following property of
these sets.

\begin{lemma}
\label{petit-osama}
If $A$ is a positively invariant set (mod $0$) of finite measure
of a negative nonsingular map $f$, then
\begin{equation}
\label{label}
\mu\left(\bigcap_{n=0}^\infty f^{-n}(A)\right)
=\mu(A).
\end{equation}
\end{lemma}

\proof
Since
$\mu(A)=\mu(A\setminus f^{-1}(A))+\mu(A\cap f^{-1}(A))$
and $A$ is positively invariant (mod $0$) one has
$\mu(A)=\mu(A\cap f^{-1}(A))$,
i.e.,
$$
\mu\left(\bigcap_{n=0}^1 f^{-n}(A)\right)=\mu(A).
$$
Now suppose that $m\in \mathbb{N}^+$ satisfies
$$
\mu\left(\bigcap_{n=0}^m f^{-n}(A)\right)=\mu(A).
$$
Since
$$
\mu\left(\bigcap_{n=0}^{m+1} f^{-n}(A)\right)=
\mu\left(\bigcap_{n=0}^m f^{-n}(A)\right)-
\mu\left(\left(\bigcap_{n=0}^m f^{-n}(A)\right)\setminus f^{-m-1}(A)\right)
$$
and
\begin{eqnarray*}
\mu\left(\left(\bigcap_{n=0}^m f^{-n}(A)\right)\setminus f^{-m-1}(A)\right)
& \leq & \mu(f^{-m}(A)\setminus f^{-m-1}(A)) \\
& = & \mu(f^{-m}(A\setminus f^{-1}(A))) \\
& = & 0
\end{eqnarray*}
because  $f$ is negative nonsingular and $A$ is positively invariant (mod $0$), one has
$\mu\left(\bigcap_{n=0}^{m+1} f^{-n}(A)\right)=\mu(A)$. Therefore
\begin{equation}
\label{macarrone}
\mu\left(\bigcap_{n=0}^m f^{-n}(A)\right)=\mu(A),
\quad\quad\forall m\in\mathbb{N},
\end{equation}
by induction.
On the other hand,
$$
\bigcap_{n=0}^\infty f^{-n}(A)=\bigcap_{m=0}^\infty
\bigcap_{n=0}^m f^{-n}(A)
$$
and
$\bigcap_{n=0}^{m+1} f^{-n}(A)\subset \bigcap_{n=0}^m f^{-n}(A)$.
As $\mu(A)<\infty$
we conclude that
$$
\mu\left(\bigcap_{n=0}^\infty f^{-n}(A)\right)
=\lim_{m\to\infty}\mu\left(\bigcap_{n=0}^{m} f^{-n}(A)\right)
\stackrel{(\ref{macarrone})}{=}
\lim_{m\to\infty}\mu(A)
=\mu(A)
$$
proving (\ref{label}).
\endproof

We use the above lemma only in the proof of
the proposition below.

\begin{prop}
\label{osama1}
Let $P$ be a measure-sensitive partition of a negative nonsingular map $f$.
Then, no $\xi\in P$ with positive finite measure is positively invariant (mod $0$).
\end{prop}

\proof
Suppose by contradiction that there is $\xi\in P$ with $0<\mu(\xi)<\infty$
which is positively invariant (mod $0$).
Taking $A=\xi$ in Lemma \ref{petit-osama}
we obtain
\begin{equation}
\label{labell}
\mu\left(\bigcap_{n=0}^\infty f^{-n}(\xi)\right)
=\mu(\xi).
\end{equation}
As $\mu(\xi)>0$ we conclude that $\bigcap_{n=0}^\infty f^{-n}(\xi)\neq\emptyset$, and so,
there is $x\in \xi$ such that $f^n(x)\in \xi$ for all $n\in\mathbb{N}$.
As $\xi\in P$ we obtain
$P(f^n(x))=\xi$ and so
$f^{-n}(P(f^n(x)))=f^{-n}(\xi)$ for all $n\in \mathbb{N}$.
Using (\ref{pn}) we get
$$
P_m(x)=\bigcap_{n=0}^mf^{-n}(\xi).
$$
Then, (\ref{pinfinity}) yields
$$
P_\infty(x)=\bigcap_{m=0}^\infty P_m(x)=\bigcap_{m=0}^\infty\bigcap_{n=0}^m f^{-n}(\xi)
=\bigcap_{n=0}^\infty f^{-n}(\xi)
$$
and so $\mu(P_\infty(x))=\mu(\xi)$ by (\ref{labell}).
Then, $\mu(\xi)=0$ by Lemma \ref{l0} since $P$ is measure-sensitive which is absurd.
This contradiction proves the result.
\endproof

We also need the following lemma
resembling a well-known property of the expansive maps.

\begin{lemma}
 \label{l2}
If $k\in \mathbb{N}^+$, then $f$ is measure-expansive if and only if $f^k$ is.
\end{lemma}

\proof
The notation
$P^f_\infty(x)$ will indicate the dependence of $P_\infty(x)$ on $f$.

First suppose that $f^k$ is an measure-expansive with measure-sensitive partition
$P$. Then, $\mu(P_\infty^{f^k}(x))=0$ for all $x\in X$ by Lemma \ref{l0}.
But by definition one has
$P_\infty^f(x)\subset P_\infty^{f^k}(x)$ so
$\mu(P_\infty^f(x))=0$ for all $x\in X$.
Therefore, $f$ is measure-expansive with measure-sensitive partition
$P$.
Conversely, suppose that $f$ is measure-expansive with expansivity constant $P$.
Consider $Q=\bigvee_{i=0}^kf^{-i}(P)$ which is a countable partition
satisfying $Q(x)=\bigcap_{i=0}^kf^{-i}(P(f^i(x)))$ by (\ref{lee1}).
Now, take $y\in Q_\infty^{f^k}(x)$.
In particular, $y\in Q(x)$ hence
$f^i(y)\in P(f^i(x))$ for every $0\leq i\leq k$.
Take $n>k$ so $n=pk+r$ for some nonnegative integers
$p$ and $0\leq r<k$.
As $y\in Q_\infty^{f^k}(x)$ one has  $f^{pk}(y)\in Q(f^{pk}(x))$
and then $f^n(y)=f^{pk+i}(y)=f^{i}(f^{pk}(y))\in P(f^i(f^{pk}(x))=P(f^n(x))$
proving $f^n(y)\in P(f^n(x))$ for all $n\in\mathbb{N}$.
Then, $y\in P_\infty(x)$ yielding
$
Q^{f^k}_\infty(x)\subset P_\infty(x).
$
Thus
$\mu(Q_\infty^{f^k}(x))=0$ for all $x\in X$ by the equivalence (ii) in Lemma \ref{l0}
since $P$ is measure-sensitive.
It follows that $f^k$ is measure-expansive with measure-sensitive partition $Q$.
\endproof

With these definitions and preliminary results we obtain the following.

\begin{theorem}
Every measure-expansive negative nonsingular map in a probability space is aperiodic*.
\end{theorem}

\proof
Suppose by contradiction that there is a measure-expansive map $f$ which is negative nonsingular
but not aperiodic*.
Then, there are a measurable set of positive measure $A$ and $n\in\mathbb{N}^+$
such that $\mu(B\setminus f^{-n}(B))=0$ for every measurable subset
$B\subset A$. It follows that every measurable subset $B\subset A$ is positively invariant (mod $0$)
with respect to $f^n$.
By Lemma \ref{l2} we can assume $n=1$.

Now, let $P$ be a measure-sensitive partition of $f$.
Clearly, since $\mu(A)>0$ there is $\xi\in P$
such that $\mu(A\cap \xi)>0$.
Taking $\eta=A\cap\xi$ we obtain that $\eta$ is positively invariant (mod $0$) with positive measure.
In addition, consider the new partition
$Q=(P\setminus \{\xi\})\cup \{\eta,\xi\setminus A\}$ which is clearly measure-sensitive
(for $P$ is). Since this partition
also carries a positively invariant (mod $0$) member of positive measure (say $\eta$)
we obtain a contradiction by
Proposition \ref{osama1}.
The proof follows.
\endproof

To finish we compare the measure-expansivity
with the notion of pairwise sensitivity in metric measure spaces introduced in p. 376 of \cite{cj}.

By a {\em metric measure space} we mean a triple $(X,d,\mu)$ where $(X,d)$ is a metric space
and $\mu$ is a measure in the corresponding Borel $\sigma$-algebra.
Hereafter the term {\em measurable} will mean {\em Borel measurable}.
The product measure in $X\times X$ will be denoted by $\mu^{\otimes 2}$.

\begin{definition}
\label{pe-measure}
A measurable map $f: X\to X$ of a metric measure space $(X,d,\mu)$ is
{\em pairwise sensitive}
if there is $\delta>0$ such that
$$
\mu^{\otimes 2}\left(\left\{(x,y)\in X\times X:\exists n\in\mathbb{N}\mbox{ such that }
d(f^n(x),f^n(y))\geq \delta\right\}\right)=1.
$$
\end{definition}

The following is a characterization of pairwise sensitivity
which is similar to one in \cite{m}.
Since this reference is not available yet we include its proof here for the sake of completeness.
By a {\em metric probability space} we mean a metric measure space
of total mass one.
Given a map $f:X\to X$ and $\delta>0$ we define the dynamical $\delta$-balls
$$
\Phi_\delta(x)=\{y\in X:d(f^n(x),f^n(y))\leq \delta,\forall n\in\mathbb{N}\},
\quad\quad\forall x\in X.
$$
\begin{lemma}
\label{th1-new}
The following properties are equivalent
for measurable maps $f$ on metric probability spaces $(X,d,\mu)$:
\begin{enumerate}
\item
$f$ is pairwise sensitive.
\item
There is $\delta>0$ such that
\begin{equation}
\label{paratodo}
\mu(\Phi_\delta(x))=0,\quad\quad
\forall x\in X.
\end{equation}
\item
There is $\delta>0$ such that
\begin{equation}
\label{quasitodo}
\mu(\Phi_\delta(x))=0,
\quad\quad
\forall \mu\mbox{-a.e. }
x\in X.
\end{equation}
\end{enumerate}
\end{lemma}

\proof
First we prove (2) and (3) are equivalent.
Indeed, we only have to prove that (3) implies (2).
Fix $\delta>0$ satisfying (\ref{quasitodo})
and suppose by contradiction that (2) fails.
Then,
there is $x_0\in X$ such that $\mu(\Phi_{\delta/2}(x_0))>0$.
Denote $X_\delta=\{x\in X:\mu(\Phi_\delta(x))=0\}$
so $\mu(X_\delta)=1$.
Since $\mu$ is a probability we obtain
$X_\delta\cap \Phi_{\frac{\delta}{2}}(x_0)\neq\emptyset$ so
there is $y_0\in \Phi_{\frac{\delta}{2}}(x_0)$ such that $\mu(\Phi_\delta(y_0))=0$.
Now if $x\in \Phi_{\frac{\delta}{2}}(x_0)$ we have $d(f^i(x),f^i(x_0))\leq \frac{\delta}{2}$
(for all $i\in \mathbb{N}$) and, since
$y_0\in \Phi_{\frac{\delta}{2}}(x_0)$, we obtain
$d(f^i(y_0),f^i(x_0))\leq \frac{\delta}{2}$ (for all $i\in \mathbb{N}$)
so
$d(f^i(x),f^i(y_0))\leq d(f^i(x),f^i(x_0))+d(f^i(x_0),f^i(y_0))\leq \frac{\delta}{2}+\frac{\delta}{2}=\delta$
(for all $i\in \mathbb{N}$)
proving $x\in \Phi_{\frac{\delta}{2}}(y_0)$.
Therefore $\Phi_{\frac{\delta}{2}}(x_0)\subset \Phi_\delta(y_0)$
so $\mu(\Phi_{\frac{\delta}{2}}(x_0))\leq \mu(\Phi_\delta(y_0))=0$ which is absurd.
This proves that (2) and (3) are equivalent.

On the other hand, given $\delta>0$ define
$\mathcal{A}_\delta$ and $f\times f:X\times X\to X\times X$ by
$$
\mathcal{A}_\delta=\{(x,y)\in X\times X:d(x,y)<\delta\}
\quad
\mbox{ and }
\quad
(f\times f)(x,y)=(f(x),f(y)).
$$
As noticed in \cite{cj}, $f$ is
pairwise sensitive if and only if there is
$\delta>0$ satisfying
\begin{equation}
\label{eq2-new}
\mu^{\otimes 2}
\left(
\bigcap_{n\in \mathbb{N}}(f\times f)^{-n}(\mathcal{A}_{\delta})
\right)=0.
\end{equation}
On the other hand, the following inequalities hold
$$
\bigcap_{n\in \mathbb{N}}(f\times f)^{-n}(\mathcal{A}_\delta)\subset
\bigcup_{x\in X}(\{x\}\times \Phi_\delta(x))
\subset
\bigcap_{n\in \mathbb{N}}(f\times f)^{-n}(\mathcal{A}_{2\delta})
$$
so
$$
F_\delta(x,y)\leq \chi_{\Phi_\delta(x)}(y)\leq F_{2\delta}(x,y),
$$
where $F_\delta$ and $\chi_C$ denotes the characteristic functions of
$\bigcap_{n\in \mathbb{N}}(f\times f)^{-n}(\mathcal{A}_\delta)$ and $C\subset X$ respectively.
Integrating the last expression we obtain
$$
\mu^{\otimes 2}
\left(
\bigcap_{n\in\mathbb{N}}(f\times f)^{-n}(\mathcal{A}_\delta)
\right)
\leq
\int_X\int_X  \chi_{\Phi_\delta(x)}(y)   d\mu(y)d\mu(x)
$$
\begin{equation}
\label{eqqq2}
\leq
\mu^{\otimes 2}
\left(
\bigcap_{n\in\mathbb{N}}(f\times f)^{-n}(\mathcal{A}_{2\delta})
\right)
\end{equation}

Now suppose that (1) holds, i.e., $f$ is pairwise sensitive.
So, there is $\delta>0$ satisfying
(\ref{eq2-new}).
Then, the second inequality in
(\ref{eqqq2}) implies $\mu(\Phi_{\frac{\delta}{2}}(x))=0$ for $\mu$-a.e. $x\in X$
whence (3) holds.
On the other hand, suppose that (2) holds, i.e.,
there is $\delta>0$ satisfying (\ref{paratodo}).
Then, the first inequality in (\ref{eqqq2}) implies (\ref{eq2-new}) so $f$ is pairwise sensitive
whence (1) holds.
This proves the result.
\endproof

By a {\em separable probability space} we mean a metric probability space
whose underlying metric space is separable.

\begin{theorem}
\label{c1}
All pairwise sensitive maps on separable probability spaces
are measure-expansive.
\end{theorem}

\proof
Let $f$ be a pairwise sensitive map
of a separable probability space $(X,d,\mu)$.
By Lemma \ref{th1-new} there is $\delta>0$ satisfying (\ref{paratodo}).
Since $(X,d)$ is separable we can select a countable covering $\{B_k:k\in I\}$
of $X$ consisting of balls of radio $\delta$,
where $I$ is either $\mathbb{N}$ or $\{0,1,\cdots, s\}$ for some $s\in \mathbb{N}$.
As usual we can transform this covering into a countable partition $P=\{\xi_k:k\in I\}$
by taking $\xi_0=B_0$ and $\xi_k=B_k\setminus \cup_{i=0}^{k-1}B_i$ for $k\geq1$.
Clearly this partition satisfies
$P_\infty(x)\subset \Phi_\delta(x)$. Then, (\ref{paratodo}) implies
$\mu(P_\infty(x))\leq \mu(\Phi_\delta(x))=0$ for every $x\in X$ so $P$ is measure-sensitive
by Lemma \ref{l0}.
\endproof

The following example shows that converse of Theorem \ref{c1} is false.

\begin{example}
\label{ex2}
An irrational circle rotation
is measure-expansive with respect to the Lebesgue measure (c.f. Example \ref{re2}
or Corollary \ref{th1})
but not pairwise sensitive with respect to that measure (c.f. \cite{cj} p. 378).
\end{example}

Recall that
a map $f: X\to X$ of a metric space
$(X,d)$ is {\em expansive}
if there is
$\delta>0$ such that
$x=y$ whenever $x,y\in X$ and $d(f^n(x),f^n(y))\leq \delta$ for all
$n\in\mathbb{N}$.

\begin{corollary}
Every measurable expansive map in a nonatomic separable probability space
is measure-expansive.
\end{corollary}

\proof
Notice that a map $f$ is expansive
if and only if there is $\delta>0$
such that $\Phi_\delta(x)=\{x\}$ for every $x\in X$.
Then, Lemma \ref{th1-new} implies that
every expansive measurable map of a nonatomic metric measure space
is pairwise sensitive.
Now apply Theorem \ref{c1}.
\endproof

\vspace{10pt}

{\small
{\em Authors' addresses}:
Instituto de Matem\'atica, Universidade Federal do Rio de Janeiro,
P. O. Box 68530, 21945-970, Rio de Janeiro, Brazil, e-mail: \texttt{morales@impa.br}.}

\end{document}